\documentclass{amsart}

\usepackage[english]{babel}
\usepackage{leftidx}
\usepackage[numbers]{natbib}
\usepackage{bigints}
\usepackage{xcolor}
\usepackage{enumitem}
\usepackage{empheq}
\usepackage{array}
\usepackage{graphicx}
\usepackage{lipsum}
\usepackage{amssymb}
\usepackage{amsthm}
\usepackage{needspace}
\usepackage[makeroom]{cancel}
\newenvironment{myproof}[2] {\paragraph{\textbf{Proof of {#1} {#2} :}}}{\hfill$\square$}

\newtheorem{theorem}{Theorem}[section]
\newtheorem{question-non}[]{}

\newtheorem{cor}[theorem]{Corollary}%[section]
\usepackage{geometry}
\title{TRIVIALITY RESULTS
FOR COMPACT $k$-YAMABE SOLITONS}

\author{Tokura, W. $^{1}$}
\address{$^{1}$ Instituto Federal Goiano, 75380-000, Av. Wilton Monteiro da Rocha, s/n, Trindade, GO, Brazil.}
\email{williamisaotokura@hotmail.com $^{1}$}

\author{Batista, E. $^{2}$}
\address{$^{2}$ Universidade Federal de Goi\'as, IME, 131, 74001-970, Goi\^ania, GO, Brazil.}
\email{edbatista@gmail.com.br $^{2}$}



\begin{document}

\begin{abstract}
In this paper, we show that any compact gradient $k$-Yamabe soliton must have constant $\sigma_{k}$-curvature. Moreover, we provide a certain condition for a compact $k$-Yamabe soliton to be gradient.
\end{abstract}
\maketitle
\section{Introduction and main results}
\label{intro}

The concept of gradient $k$-Yamabe soliton, introduced in the celebrated work \cite{catino2012global}, corresponds to a natural generalization of gradient Yamabe solitons. We recall that a Riemannian manifold $(M^n, g)$ is a \textit{$k$-Yamabe soliton} if it admits a constant $\lambda\in\mathbb{R}$ and a vector field $X\in \mathfrak{X}(M)$ satisfying the equation 
\begin{equation}\label{def1}
    \frac{1}{2}\mathcal{L}_{X}g=2(n-1)(\sigma_{k}-\lambda)g,
\end{equation}
where $\mathcal{L}_{X}g$ and $\sigma_{k}$ stand, respectively, for the Lie derivative of $g$ in the direction of $X$ and the $\sigma_{k}$-curvature of $g$.  Recall that, if we denote by $\lambda_{1},\lambda_{2},\dots,\lambda_{n}$ the
eigenvalues of the symmetric endomorphism $g^{-1}A_{g}$, where $A_{g}$ is the Schouten tensor defined by
\begin{equation*}
A_{g}=\frac{1}{n-2}\left(Ric_{g}-\frac{scal_{g}}{2(n-1)}g\right),
\end{equation*}
then the $\sigma_{k}$-curvature of $g$ is defined as the $k$-th symmetric elementary function of $\lambda_{1},\dots,\lambda_{n}$, namely
\begin{equation*}
\sigma_{k}=\sigma_{k}(g^{-1}A_{g})=\sum_{i_{1}<\dots i_{k}}\lambda_{i_{1}}\cdot \dots\cdot \lambda_{i_{k}}, \quad \text{for}\quad 1\leq k\leq n.
\end{equation*}
Since $\sigma_{1}$ is the trace of $g^{-1}A_{g}$, the $1$-Yamabe solitons simply correspond to gradient Yamabe solitons \cite{chow1992yamabe,daskalopoulos2013classification, di2008yamabe, hamilton1988ricci, ma2012remarks,tokura2018warped}. For simplicity, the soliton will be denoted by $(M^{n}, g, X, \lambda)$. It may happen that $X=\nabla f$ is the gradient field of a smooth real function $f$ on $M$, in which case the soliton $(M^{n}, g, \nabla f, \lambda)$ is referred to as a \textit{gradient $k$-Yamabe soliton}. Equation \eqref{def1} then becomes
\begin{equation}\label{eq fundamental}
\nabla^2 f=2(n-1)(\sigma_{k}-\lambda)g,
\end{equation}
where $\nabla^2 f$ is the Hessian of $f$. Moreover, when either $f$ is a constant function or $X$ is a Killing vector field, the soliton is called \textit{trivial} and, in this case, the metric $g$ is of constant $k$-curvature $\sigma_{k}=\lambda$.

In recent years, much efforts have been devoted to study the geometry of $k$-Yamabe solitons. For instance, Hsu in \cite{hsu2012note} shown that any compact gradient $1$-Yamabe soliton is trivial. For $k>1$, the extension of the previous result was investigated by Catino et al. \cite{catino2012global}, and Bo et al. \cite{bo2018k}. In  \cite{catino2012global}, the authors proved that any compact, gradient $k$-Yamabe soliton with nonnegative Ricci tensor is trivial. On the other hand, the authors in \cite{bo2018k} showed that any compact, gradient $k$-Yamabe soliton with constant negative scalar curvature must be trivial.

In this paper, we extend the above results as follows.
%The gradient soliton is called \textit{trivial} if $h$ is a constant function.

\begin{theorem}\label{T1}Any compact gradient $k$-Yamabe soliton $(M^n,g, \nabla f,\lambda)$ is trivial, i.e., has constant $\sigma_{k}$-curvature $\sigma_{k}=\lambda$.
\end{theorem}

In the scope of $k$-Yamabe solitons, we provide the following extension of Theorem 1.3 in \cite{bo2018k}.
%From the analogous arguments present in the demonstration of previous Theorem, we derive the following result for $k$-Yamabe solitons.

\begin{theorem}\label{T4}The compact $k$-Yamabe soliton $(M^n,g,X,\lambda)$ is trivial if one of the following conditions holds:
\begin{itemize}
    \item [\textup{(a)}] $k=1$.
    \vspace{0,1cm}
    \item [\textup{(b)}] $k\geq2$ and $(M^n,g)$ is locally conformally flat.
\end{itemize}
\end{theorem}

%Since $k$-Yamabe solitons with $k=1$ are trivial, our next triviallity results are given in the case $k\geq2$.

The Hodge-de Rham decomposition theorem (see \cite{aquino2011some, warner2013foundations}),  shows that any vector field $X$ on a compact oriented Riemannian manifold $M$ can be decompose as follows: 
\begin{equation}\label{Hodge} X = \nabla h + Y,
\end{equation}
where $h$ is a smooth function on $M$ and $Y\in \mathfrak{X}(M)$ is a free divergence vector field. Indeed, just consider the $1$-form $X^{\flat}$. Hence applying the Hodge-de Rham theorem, we decompose $X^{\flat}$ as follows:
\[X^{\flat}=d\alpha+\delta\beta+\gamma.\]
Taking $Y = (\delta\beta+\gamma)^{\sharp}$ and $(d\alpha)^{\sharp}=\nabla h$ we arrive at the desired result.

Now we notice that the same result obtained in \cite{pirhadi2017almost} for compact almost Yamabe solitons also works for compact $k$-Yamabe solitons. More precisely, we have the following theorem.

\begin{theorem}\label{T2}The compact $k$-Yamabe soliton $(M^n,g,X,\lambda)$ is gradient if, and only if,
\[\int_{M^n}Ric(\nabla h,Y)dv_{g} \leq0,\]
where $h$ and $Y$ are the Hodge-de Rham decomposition components of $X$.
\end{theorem}

As a consequence of Theorem \ref{T1} and Theorem \ref{T2}, we derive the following triviality result.
%\begin{cor}Any compact $k$-Yamabe soliton $(M^n,g,X,\lambda)$ satisfying \[\int_{M^n}\textup{Ric}(\nabla h,Y)\leq0,\]
 %where $h$ and $Y$ are the Hodge–de Rham decomposition components of $X$ has constant $\sigma_{k}$-curvature $\sigma_{k}=\lambda$.
%\end{cor}

\begin{cor}\label{corr}Let $(M^n,g,X,\lambda)$ be a compact $k$-Yamabe soliton $(k\ge2)$ and $X=\nabla h+Y$ the Hodge-de Rham decomposition of $X$. If 
\[\int_{M^n}Ric(\nabla h,Y)dv_{g}\leq0,\]
then $(M^n, g)$ is a trivial $k$-Yamabe soliton.
\end{cor}

An immediate consequence of the above corollary is the next result.

\begin{cor}\label{cor}Any compact $k$-Yamabe soliton $(M^n,g,X,\lambda)$ with $k\ge2$ and nonpositive Ricci curvature is trivial.
\end{cor}

Finally, taking into account the $L^{2}(M)$ orthogonality of the Hodge-de Rham decomposition, we obtain.

\begin{theorem}\label{1.6}Let $(M^n,g,X,\lambda)$ be a compact $k$-Yamabe soliton $(k\ge2)$ and   $X=\nabla h+Y$ the Hodge-de Rham decomposition of $X$. If 
\[\int_{M^n}g(\nabla h,X)dv_{g}\leq0,\]
then $(M^n, g)$ is a trivial $k$-Yamabe soliton.

\end{theorem}

%\begin{theorem}Let $(M^n,g,X,\lambda)$ be a compact $k$-Yamabe soliton $(k\ge2)$ and   $X=\nabla h+Y$ the Hodge–de Rham decomposition of $X$. If 
%\[\int_{M^n}g(\nabla h,X)dv_{g}\leq0.\]
%Then $(M^n, g)$ is a trivial $k$-Yamabe soliton.

%\end{theorem}

%\begin{theorem}Let $(M^n,g,X,\lambda)$ be a complete noncompact $k$-Yamabe soliton $(k\ge2)$ with constant scalar curvature. If $|X|\in L^{1}(M)$
%Then $(M^n, g)$ is a trivial $k$-Yamabe soliton.
%\end{theorem}

\section{Proofs}

\begin{myproof}{Theorem}{\ref{T4}} If $k=1$, then $(M^n,g)$ is a Yamabe soliton and the result is well known from \cite{di2008yamabe}. Now, consider $k\geq2$ and suppose $(M^n,g)$ locally conformally flat. It was proved in \cite{han2006kazdan, viaclovsky2000some} that, on a compact, locally conformally flat, Riemannian manifold,
one has 
\[\int_{M^n}g(\nabla \sigma_{k},X)dv_{g}=0,\]
for every conformal Killing vector field $X$ on $(M^n,g)$. From the structure equation \eqref{def1}, we know that $X$ is a conformal Killing vector field; hence, it follows that
\begin{equation}\label{9090}
    0=\int_{M^n}g(\nabla \sigma_{k},X)dv_{g}=-\int_{M^n}\sigma_{k}(div X)dv_{g}=-2n(n-1)\int_{M^n}\sigma_{k}(\sigma_{k}-\lambda)dv_{g},
\end{equation}
where in the second equality we have used the divergence theorem. On the other hand, again from the divergence theorem, we obtain
\begin{equation}\label{8080}0=\int_{M^n}div X dv_{g}=2n(n-1)\int_{M^n}(\sigma_{k}-\lambda)dv_{g}.
\end{equation}
Jointly equations \eqref{9090} and \eqref{8080}, we conclude that
\[2n(n-1)\int_{M^n}(\sigma_{k}-\lambda)^2 dv_{g}=0,\]
which implies that $\sigma_{k}=\lambda$ and $\mathcal{L}_{X}g=0$. Hence $(M^n,g)$ is trivial.

\end{myproof}

\begin{myproof}{Theorem}{\ref{T1}} If $k=1$, then $(M^n,g)$ is a gradient Yamabe soliton and the result is well known from \cite{hsu2012note}. Now, consider $k\geq2$ and suppose by contradiction that $f$ is nonconstant. From Theorem 1.1 of \cite{catino2012global}, we obtain that $(M^n,g)$ is rotationally symmetric and $M^{n}\setminus \{N,S\}$ is locally conformally flat. Here $N,S$  corresponds to the extremal points of  $f$ in $M$. From the structure equation \eqref{eq fundamental}, we know that $\nabla f$ is a conformal Killing vector field; hence, we can apply Theorem 5.2 of \cite{viaclovsky2000some} to deduce
\begin{equation}\label{1221}
    0=\int_{M^{n}\setminus\{N,S\}}g(\nabla\sigma_{k},\nabla f)dv_{g}=\int_{M^n}g(\nabla\sigma_{k},\nabla f)dv_{g}=-2n(n-1)\int_{M^n}\sigma_{k}(\sigma_{k}-\lambda)dv_{g},
\end{equation}
where in the last equality we have used the divergence theorem. On the other hand, again from the divergence theorem, we get
\begin{equation}\label{808080}0=\int_{M^n}\Delta f dv_{g}=2n(n-1)\int_{M^n}(\sigma_{k}-\lambda)dv_{g}.
\end{equation}
Jointly equations \eqref{1221} and \eqref{808080}, we conclude that
\[2n(n-1)\int_{M^n}(\sigma_{k}-\lambda)^2 dv_{g}=0,\]
which implies that $\sigma_{k}=\lambda$ and $f$ is harmonic. Since $M^n$ is compact, $f$ is a constant, which leads to a contradiction. This proves that $f$ is constant.

\end{myproof}

\begin{myproof}{Theorem}{\ref{T2}}
From the Hodge-de Rham decomposition \eqref{Hodge}, we deduce that
\begin{equation}\label{t11111}\frac{1}{2}\mathcal{L}_{Y}g=\frac{1}{2}\mathcal{L}_{X}g-\frac{1}{2}\mathcal{L}_{\nabla h}g=2(n-1)(\sigma_{k}-\lambda)g-\nabla^{2}h.
\end{equation}
Therefore, to prove that $(M^n,g)$ admits a gradient $k$-Yamabe soliton structure, it is necessary and sufficient to show that $\mathcal{L}_{Y}g=0$. From \eqref{t11111}, we arrive that
\begin{equation}\label{t12}
    \begin{split}
        \frac{1}{4}\int_{M^n}|\mathcal{L}_{Y}g|^2 dv_{g}&=\int_{M^n}\left[4n(n-1)^2(\sigma_{k}-\lambda)^2-4(n-1)g\left(\nabla^2h,(\sigma_{k}-\lambda)g\right)+|\nabla^2 h|^2 \right]dv_{g}\\
        &=\int_{M^n}\left[|\nabla^2 h|^2-4n(n-1)^2(\sigma_{k}-\lambda)^2\right]dv_{g}.
    \end{split}
\end{equation}
We are going to compute the right-hand side of \eqref{t12} using the following identity %order to compute the right side of \eqref{t12}, let us first calculate the following integral
\begin{equation}\label{tt1}\int_{M^n}2Ric(\nabla h,Y)dv_{g}=\int_{M^n}\left[Ric(X,X)-Ric(\nabla h,\nabla h)-Ric(Y,Y)\right]dv_{g}.
\end{equation}

Taking the divergence of \eqref{t11111}, we get
\begin{equation}\label{t2}
\begin{split}
\frac{1}{2}div(\mathcal{L}_{Y}g)(Y)&=\frac{1}{2}div(\mathcal{L}_{X}g)(Y)-\frac{1}{2}div(\mathcal{L}_{\nabla h}g)(Y)\\
&=2(n-1)div(\sigma_{k}-\lambda)(Y)-\frac{1}{2}div(\mathcal{L}_{\nabla h}g)(Y)\\
&=2(n-1)g(\nabla\sigma_{k},Y)-\frac{1}{2}div(\mathcal{L}_{\nabla h}g)(Y).
\end{split}
\end{equation}
Hence, from the Bochner formula (see Lemma 2.1 of \cite{petersen2009rigidity}), we can express \eqref{t2} as follows
\begin{equation}\label{t4}
    \frac{1}{2}\Delta|Y|^2-|\nabla Y|^2+Ric(Y,Y)=4(n-1)g(\nabla\sigma_{k},Y)-2Ric(\nabla h, Y)-2g(\nabla \Delta h,Y),
\end{equation}
and using the compactness of $M^n$, we arrive at
equation
\begin{equation}\label{dd}
    \int_{M^n}2Ric(\nabla h,Y)dv_{g}=\int_{M^n}\left[|\nabla Y|^2-Ric(Y,Y)\right] dv_{g}.
\end{equation}
On the other hand, the same argument as above shows that
\begin{equation}\label{11}
    \frac{1}{2}\Delta|X|^2-|\nabla X|^2+Ric(X,X)=-2(n-1)(n-2)g(\nabla \sigma_{k},X).
\end{equation}
Since
\begin{equation*}
    \begin{split}
        \int_{M^n}|\nabla X|^2dv_{g}&=\int_{M^n}\big{[}|\nabla^2 h|^2+|\nabla Y|^2+2g(\nabla\nabla h,\nabla Y)\big{]}dv_{g}\\
        &=\int_{M^n}\big{[}|\nabla^2 h|^2+|\nabla Y|^2-2g(\nabla\Delta h+Ric(\nabla h), Y)\big{]}dv_{g}\\
        &=\int_{M^n}\big{[}|\nabla^2 h|^2+|\nabla Y|^2-2Ric(\nabla h, Y)\big{]}dv_{g},\\
    \end{split}
\end{equation*}
we may integrate \eqref{11} over $M^n$ to deduce
\begin{equation}\label{dd2}
    \begin{split}
        \int_{M^n}Ric(X,X)dv_{g}&=\int_{M^n}\left[|\nabla X|^2-2(n-1)(n-2)g(\nabla \sigma_{k},X)\right]dv_{g}\\
        &=\int_{M^n}\big{[}|\nabla^2 h|^2+|\nabla Y|^2-2Ric(\nabla h, Y)+4n(n-1)^2\times\\
        &\qquad\times(n-2)(\sigma_{k}-\lambda)^2\big{]}dv_{g}.\\
    \end{split}
\end{equation}
Again, the same argument based on Lemma 2.1 of \cite{petersen2009rigidity}, allow us to deduce that
\begin{equation}\label{dd3}
\int_{M^n}Ric(\nabla h,\nabla h)dv_{g}=\int_{M^n}\left[4n^2(n-1)^2(\sigma_{k}-\lambda)-|\nabla^2 h|^2\right]dv_{g}.
\end{equation}

Now, replacing back \eqref{dd}, \eqref{dd2} and \eqref{dd3} into \eqref{tt1}, we get
\[\int_{M^n}\left[|\nabla^2 h|-4n(n-1)^2(\sigma_{k}-\lambda)^2\right]dv_{g}=\int_{M^n}Ric(\nabla h,Y)dv_{g},\]
which combining with \eqref{t12} produce the desired result.

\end{myproof}

\begin{myproof}{Theorem}{\ref{1.6}}Since the Hodge-de Rham decomposition is
orthogonal on $L^{2}(M)$, we get
\begin{equation*}
    \int_{M^n}g(\nabla h, X)dv_{g}=\int_{M^n}g(\nabla h, \nabla h+Y)dv_{g}=\int_{M^n}|\nabla h|^{2}dv_{g}.
\end{equation*}
Therefore, if
\begin{equation*}
    \int_{M^n}g(\nabla h, X)dv_{g}\leq0,
\end{equation*}
we obtain that $\nabla h=0$  and, consequently, $X=Y$. Now, since $Y$ is a free divergence vector field, we deduce
\[0=div Y=div X=2n(n-1)(\sigma_{k}-\lambda),\]
which implies that $\sigma_{k}=\lambda$ and $\mathcal{L}_{X}g=0$, hence, trivial.

\end{myproof}

%\begin{cor}Let $(M^n,g,X,\lambda)$ be a compact $k$-Yamabe soliton satisfying \[\int_{M^n}\textup{Ric}(\nabla f,Y)\leq0,\]
%where $f$ and $Y$ are the Hodge–de Rham decomposition components of $X$. Then $X$ is a Killing vector field.
%\end{cor}


\begin{thebibliography}{10}
	
	\bibitem{aquino2011some}
	C.~Aquino, A.~Barros, and E.~Ribeiro.
	\newblock Some applications of the hodge-de rham decomposition to ricci
	solitons.
	\newblock {\em Results in Mathematics}, 60(1-4):245, 2011.
	
	\bibitem{bo2018k}
	L.~Bo, P.~T. Ho, and W.~Sheng.
	\newblock The k-yamabe solitons and the quotient yamabe solitons.
	\newblock {\em Nonlinear Analysis}, 166:181--195, 2018.
	
	\bibitem{catino2012global}
	G.~Catino, C.~Mantegazza, and L.~Mazzieri.
	\newblock On the global structure of conformal gradient solitons with
	nonnegative ricci tensor.
	\newblock {\em Communications in Contemporary Mathematics}, 14(06):1250045,
	2012.
	
	\bibitem{chow1992yamabe}
	B.~Chow.
	\newblock The yamabe flow on locally conformally flat manifolds with positive
	ricci curvature.
	\newblock {\em Communications on pure and applied mathematics},
	45(8):1003--1014, 1992.
	
	\bibitem{daskalopoulos2013classification}
	P.~Daskalopoulos and N.~Sesum.
	\newblock The classification of locally conformally flat yamabe solitons.
	\newblock {\em Advances in Mathematics}, 240:346--369, 2013.
	
	\bibitem{di2008yamabe}
	L.~F. Di~Cerbo and M.~M. Disconzi.
	\newblock Yamabe solitons, determinant of the laplacian and the uniformization
	theorem for riemann surfaces.
	\newblock {\em Letters in Mathematical Physics}, 83(1):13--18, 2008.
	
	\bibitem{hamilton1988ricci}
	R.~S. Hamilton.
	\newblock The ricci flow on surfaces, mathematics and general relativity (santa
	cruz, ca, 1986), 237--262.
	\newblock {\em Contemp. Math}, 71:301--307, 1988.
	
	\bibitem{han2006kazdan}
	Z.-C. Han.
	\newblock A kazdan--warner type identity for the $\sigma_k$ curvature.
	\newblock {\em Comptes Rendus Mathematique}, 342(7):475--478, 2006.
	
	\bibitem{hsu2012note}
	S.-Y. Hsu.
	\newblock A note on compact gradient yamabe solitons.
	\newblock {\em Journal of Mathematical Analysis and Applications},
	388(2):725--726, 2012.
	
	\bibitem{ma2012remarks}
	L.~Ma and V.~Miquel.
	\newblock Remarks on scalar curvature of yamabe solitons.
	\newblock {\em Annals of Global Analysis and Geometry}, 42(2):195--205, 2012.
	
	\bibitem{petersen2009rigidity}
	P.~Petersen and W.~Wylie.
	\newblock Rigidity of gradient ricci solitons.
	\newblock {\em Pacific journal of mathematics}, 241(2):329--345, 2009.
	
	\bibitem{pirhadi2017almost}
	V.~Pirhadi and A.~Razavi.
	\newblock On the almost quasi-yamabe solitons.
	\newblock {\em International Journal of Geometric Methods in Modern Physics},
	14(11):1750161, 2017.
	
	\bibitem{tokura2018warped}
	W.~Tokura, L.~Adriano, R.~Pina, and M.~Barboza.
	\newblock On warped product gradient yamabe solitons.
	\newblock {\em Journal of Mathematical Analysis and Applications}, 2018.
	
	\bibitem{viaclovsky2000some}
	J.~A. Viaclovsky.
	\newblock Some fully nonlinear equations in conformal geometry.
	\newblock {\em AMS IP studies in advanced mathematics}, 16:425--434, 2000.
	
	\bibitem{warner2013foundations}
	F.~W. Warner.
	\newblock {\em Foundations of differentiable manifolds and Lie groups},
	volume~94.
	\newblock Springer Science \& Business Media, 2013.
	
\end{thebibliography}
\end{document}